\documentclass[12pt, reqno]{amsart}
 \usepackage{amsmath, amsthm, amscd, amsfonts, amssymb, graphicx, color, float}
\usepackage[bookmarksnumbered, colorlinks, plainpages]{hyperref}

\newtheorem{theorem}{Theorem}[section]
\newtheorem{lemma}[theorem]{Lemma}

\newtheorem{corollary}[theorem]{Corollary}
\theoremstyle{definition}

\theoremstyle{remark}

\numberwithin{equation}{section}

\begin{document}

\title{Inequalities for trace on $\tau$-measurable operators }
\author[M.S. Moslehian, Gh. Sadeghi]{Mohammad Sal Moslehian$^1$ and Ghadir Sadeghi$^2$}

\address{$^1$ Department of Pure Mathematics, Center of Excellence in
Analysis on Algebraic Structures (CEAAS), Ferdowsi University of
Mashhad, P. O. Box 1159, Mashhad 91775, Iran}
\email{moslehian@um.ac.ir,
moslehian@member.ams.org}

\address{$^2$ Department of Mathematics and Computer
Sciences, Hakim Sabzevari University, P.O. Box 397, Sabzevar, Iran.}
\email{ghadir54@gmail.com, g.sadeghi@hsu.ac.ir}

%======================================================================

\begin{abstract}
Let $\mathfrak{M}$ be a semifinite von Neumann algebra on a Hilbert
space equipped with a faithful normal semifinite trace $\tau$.  A
closed densely defined operator $x$ affiliated with $\mathfrak{M}$
is called $\tau$-measurable if there exists a number $\lambda \geq
0$ such that $\tau \left(e^{|x|}(\lambda,\infty)\right)<\infty$. A
number of useful inequalities, which are known for the trace on
Hilbert space operators, are extended to trace on $\tau$-measurable
operators. In particular, these inequalities imply Clarkson
inequalities for $n$-tuples of $\tau$-measurable operators. A
general parallelogram law for $\tau$-measurable operators are given
as well.
\end{abstract}

\maketitle

%======================================================================
\section{Introduction and preliminaries}
Let $\mathfrak{M}$ be a semifinite von Neumann algebra on a Hilbert
space $\mathfrak{H}$, with unit element $\textbf{1}$, equipped with
a faithful normal semifinite trace $\tau$. For standard facts
concerning von Neumann algebras, we refer the reader to \cite{D}.  A
closed densely defined linear operator $x:\mathcal{D}(x)\subseteq
\mathfrak{H\to \mathfrak{H}}$ is called affiliated with
$\mathfrak{M}$ if $ux=xu$ for all unitaries $u$ in the commutant
$\mathfrak{M}^{'}$ of $\mathfrak{M}$. If $x$ is in the algebra
$\mathcal{B}(\mathfrak{H})$ of all bounded linear operators on the
Hilbert space $\mathfrak{H}$, then $x$ is affiliated with
$\mathfrak{M}$ if and only if $x \in \mathfrak{M}$. A closed densely
defined operator $x$ affiliated with $\mathfrak{M}$ is called
$\tau$-measurable if there exists a number $\lambda \geq 0$ such
that
\begin{eqnarray*}
\tau \left(e^{|x|}(\lambda,\infty)\right)<\infty\,,
\end{eqnarray*}
where $|x|=(x^*x)^{1/2}$ and $e^a$ denotes the spectral measure of
the self-adjoint operator $a$ which is a $\sigma$-additive (w.r.t.
the strong operator topology) from the Borel $\sigma$-algebra of
$\mathbb{R}$ into the orthogonal projections. The collection of all
$\tau$-measurable operators is denoted by
$\widetilde{\mathfrak{M}}$. With the sum and product defined as the
respective closure of the algebraic sum and product, it is well
known that  $\widetilde{\mathfrak{M}}$ is a $*$-algebra with respect
to the operations of strong sum and strong product and taking
adjoint; see \cite{N}. For $x\in\widetilde{\mathfrak{M}}$, the
generalized singular value function $\mu(x): [0,\infty] \to
[0,\infty]$ is defined by
\begin{eqnarray*}
\mu_t(x)=\inf\left\{\lambda\geq0 :\quad \tau \left(e^{|x|}(\lambda,\infty)\right)\leq t\right\}, \quad t\geq 0.
\end{eqnarray*}
The generalized singular value function $\mu(x)$ is decreasing
right-continuous. Moreover,
\begin{eqnarray*}
\mu(uxv)\leq \|u\| \|v\| \mu(x)
\end{eqnarray*}
for all $u,v \in \mathfrak{M}$ and $x\in \widetilde{\mathfrak{M}}$.
Moreover,
\begin{eqnarray*}
\mu(f(x))=f(\mu(x))
\end{eqnarray*}
whenever $0\leq x\in \widetilde{\mathfrak{M}}$ and $f$ is an
increasing continuous function on $[0,\infty)$ satisfying $f(0)=0$.
The space $\widetilde{\mathfrak{M}}$ is a partially ordered vector
space under the ordering  $x\geq0$ defined by $\langle x\xi, \xi
\rangle \geq 0$, $\xi \in \mathcal{D}(x)$. The trace $\tau$ on
$\mathfrak{M}^+$ extends uniquely to an additive, positively
homogeneous, unitarily invariant and normal functional
$\widetilde{\tau}:\widetilde{\mathfrak{M}}\to[0, \infty]$, which is
given by $\widetilde{\tau}(x)=\int_0^{\infty}\mu_t(x)dt$,
$x\in\mathfrak{M}^+$ \cite{DDP3}. This extension is also denoted by
$\tau$. Further,
\begin{eqnarray}\label{EQ}
\tau(f(x))=\int_{0}^{\infty}f(\mu_t(x))dt
\end{eqnarray}
whenever $0\leq x\in \widetilde{\mathfrak{M}}$ and $f$ is
non-negative Borel function which is bounded on a neighborhood of
$0$ and satisfies $f(0)=0$.
%---------------------------------------------------------------------------------------------------------------------------------------------------------------%
For $0<p<\infty$, $L^p(\mathfrak{M},\tau)$ is defined as the set of all densely defined closed operators $x$ affiliated with $\mathfrak{M}$ such that
$$\|x\|_p=\tau(|x|^p)^{\frac{1}{p}}=\left(\int_0^{\infty}\mu_t(x)^pdt\right)^{\frac{1}{p}}<\infty.$$
For further details and proofs, we refer the reader to \cite{Y,FK,
DDP,P, SAD}. Important special cases of these noncommutative spaces
are usual $L^p$-spaces and the Schatten $p$-classes $\mathcal{C}_p$.

The classical Clarkson inequalities, for schatten $p$-norms of
Hilbert space operators, assert that if
$A,B\in\mathcal{B}(\mathfrak{H})$, then
\begin{eqnarray}\label{ca1}
2(\|A\|^p_p+\|B\|^p_p)\leq\|A+B\|^p_p+\|A-B\|^p_p\leq2^{p-1}(\|A\|^p_p+\|B\|^p_p)
\end{eqnarray}
for $2\leq p<\infty$, and
\begin{eqnarray}\label{ca2}
2^{p-1}(\|A\|^p_p+\|B\|^p_p)\leq\|A+B\|^p_p+\|A-B\|^p_p\leq2(\|A\|^p_p+\|B\|^p_p)
\end{eqnarray}
for $0<p\leq2$ (see \cite{HK1}). Bahatia and Kittaneh \cite{BK1}
generalized the Clarkson inequalities (\ref{ca1}) and (\ref{ca2})
for $n$-tuples operators as follows. If $A_0,\ldots,
A_{n-1}\in\mathcal{B}(\mathfrak{H})$ and
$\omega_0,\ldots,\omega_{n-1}$ are $n$ roots of unity with
$\omega_j=e^{2\pi ij/n}$, $0\leq j\leq n-1$, then
\begin{eqnarray}\label{ca3}
n\sum_{j=0}^{n-1}\|A_j\|_p^p\leq\sum_{k=0}^{n-1}\left\|\sum_{j=0}^{n-1}\omega_j^kA_j\right\|_p^p\leq n^{p-1}\sum_{j=0}^{n-1}\|A_j\|_p^p
\end{eqnarray}
for $2\leq p\leq\infty$, and
\begin{eqnarray}\label{ca4}
 n^{p-1}\sum_{j=0}^{n-1}\|A_j\|_p^p\leq\sum_{k=0}^{n-1}\left\|\sum_{j=0}^{n-1}\omega_j^kA_j\right\|_p^p\leq n\sum_{j=0}^{n-1}\|A_j\|_p^p
\end{eqnarray}
for $0<p<\infty$. Related Clarkson inequalities for $n$-tuples of
operators have been recently given by Kissin \cite{K} and Hirzalleh
and Kitaneh \cite{HK2}. These inequalities have been found to be
very powerful tools in operator theory and in mathematical physics
(see, e.g., \cite{S}). In \cite{FK}, Fack and Kosaki proved Clarkson
inequalities for measurable operators in Haagerup $L^p$-spaces. Let
$x$ be a $\tau$-measurable operator in $L^p(\mathfrak{M},\tau)$.
Since $\|x\|^p_p=\tau(|x|^p)$, for $0<p<\infty$, our generalization
of the inequalities (\ref{ca1}) and (\ref{ca2}) will be clear if we
rewrite them as
\begin{eqnarray}\label{ca15}
2[\tau(|x|^p)+\tau(|y|^p)]\leq\tau(|x+y|^p)+\tau(|x-y|^p)\leq2^{p-1}[\tau(|x|^p)+\tau(|y|^p)]
\end{eqnarray}
for $2\leq p<\infty$, and
\begin{eqnarray}\label{ca16}
2^{p-1}[\tau(|x|^p)+\tau(|y|^p)]\leq\tau(|x+y|^p)+\tau(|x-y|^p)\leq2[\tau(|x|^p)+\tau(|y|^p)]
\end{eqnarray}
for $0< p\leq2$.\\
There are several extensions of the classical parallelogram law in
the literature. Generalizations of the parallelogram law for the
Schatten $p$-norms have been given in the form of the celebrated
Clarkson inequalities; see \cite{HK2, MO, MST} and references cited
therein.

A number of useful inequalities relating the traces of operators on
a Hilbert space are known when the trace is defined in the usual
way. In this paper, we present some inequalities for trace on
$\tau$-measurable operators. These inequalities to trace generalized
Clarkson inequalities for $n$-tuples of $\tau$-measurable operators.
In fact, we give natural generalizations of inequalities
(\ref{ca15}) and (\ref{ca16}) for $n$-tuples of $\tau$-measurable
operators and show that the power functions $\varphi(t)=t^p$ can be
replaced by more general classes of functions. A general
parallelogram law for $\tau$-measurable operators are given as well.

%-----------------------------------------------------------------------------------------------------------------------------------------------------------------%
\section{Trace inequalities for $n$-tuple of
$\tau$-measurable operators}

In this section, we give new trace inequalities for $n$-tuple of
$\tau$-measurable operators that considerably generalize
inequalities (\ref{ca15}) and (\ref{ca16}). Before we give the main
theorems of this section, we need the following lemma, which clearly
is an extension of \cite[proposition 4.6]{FK}) for $n$-tuples of
$\tau$-measurable operators.
%---------------------------------------------------------------------------------------------------------------------------------------------------------%
\begin{lemma}\label{FK}
Let $x_0,\ldots, x_{n-1}$ be positive $\tau$-measurable operators
and $\alpha_0,\ldots, \alpha_{n-1}$ be positive real numbers such
that $\sum^{n-1}_{j=0}\alpha_j=1$. Then, for every continuous
increasing function $f$
on $\mathbb{R}_+$ with $f(0)=0$, \\
\begin{eqnarray}\label{FK1}
\tau\left(f\left(\sum_{j=0}^{n-1}\alpha_jx_j\right)\right)\leq\sum_{j=0}^{n-1}\alpha_j\tau(f(x_j)) \hspace{.3cm} \mbox{when $f$ is convex}.
\end{eqnarray}
\begin{eqnarray}\label{FK2}
\sum_{j=0}^{n-1}\alpha_j\tau(f(x_j))\leq\tau\left(f\left(\sum_{j=0}^{n-1}\alpha_jx_j\right)\right)\hspace{.3cm} \mbox{when $f$ is concave}.
\end{eqnarray}
\begin{eqnarray}\label{FK3}
\sum_{j=0}^{n-1}\tau\left(f\left(x_j\right)\right)\leq\tau \left(f\left(\sum_{j=0}^{n-1}x_j\right)\right) \hspace{.3cm} \mbox{when $f$ is convex}.
\end{eqnarray}
\begin{eqnarray}\label{FK4}
\tau \left(f\left(\sum_{j=0}^{n-1}x_j\right)\right)\leq\sum_{j=0}^{n-1}\tau\left(f\left(x_j\right)\right) \hspace{.3cm} \mbox{when $f$ is concave}.
\end{eqnarray}
\end{lemma}
%--------------------------------------------------------------------------------------------------------------------------------------------------------------%
Now we are in a position to present our main result, providing a
generalization of inequality (\ref{ca15}) for $n$-tuples of
$\tau$-measurable operators.
%-----------------------------------------------------------------------------------------------------------------------------------------------------------%
\begin{theorem}\label{MT1}
Let $x_0,\ldots,x_{n-1}$ be $\tau$-measurable operators and
$\alpha_0,\ldots,\alpha_{n-1}$ be positive real numbers such that
$\Sigma_{j=0}^{n-1}\alpha_j=1$. Let $\varphi$ be a continuous
increasing function on $\mathbb{R}_+$ with $\varphi(0)=0$ such that
$\psi(t)=\varphi(\sqrt{t})$ is convex. Then
\begin{eqnarray}\label{ca5}
\tau\left(\varphi\left(\left|\sum_{j=0}^{n-1}\alpha_jx_j\right|\right)\right)+
\sum_{0\leq j<k\leq
n-1}\tau\left(\varphi\left(\sqrt{\alpha_j\alpha_k}\left|x_j-x_k\right|\right)\right)\leq\sum_{j=0}^{n-1}\alpha_j\tau(\varphi(|x_j|)).
\end{eqnarray}
\begin{proof}
The following useful identity can be easily verified,
\begin{eqnarray}\label{ID1}
\left|\sum_{j=0}^{n-1}\alpha_jx_j\right|^2+\sum_{0\leq j<k\leq n-1}\alpha_j\alpha_k|x_j-x_k|^2=\sum_{j=0}^{n-1}\alpha_j|x_j|^2.
\end{eqnarray}
Now, by using identity (\ref{ID1}) and inequalities (\ref{FK1}),
(\ref{FK3}) we get
\begin{eqnarray*}
\sum_{j=0}^{n-1}\alpha_j\tau(\varphi(|x_j|))&=&\sum_{j=0}^{n-1}\alpha_j\tau(\psi(|x_j|^2))\\
&\geq&\tau\left(\psi\left(\sum_{j=0}^{n-1}\alpha_j|x_j|^2\right)\right)\\
&=&\tau\left(\psi\left(\left|\sum_{j=0}^{n-1}\alpha_jx_j\right|^2+\sum_{0\leq j<k\leq n-1}\alpha_j\alpha_k|x_j-x_k|^2\right)\right)\\
&\geq&\tau\left(\psi\left(\left|\sum_{j=0}^{n-1}\alpha_jx_j\right|^2\right)\right)+
\sum_{0\leq j<k\leq n-1}\tau\left(\psi\left(\alpha_j\alpha_k|x_j-x_k|^2\right)\right)\\
&=&\tau\left(\varphi\left(\left|\sum_{j=0}^{n-1}\alpha_jx_j\right|\right)\right)+
\sum_{0\leq j<k\leq n-1}\tau\left(\varphi\left(\sqrt{\alpha_j\alpha_k}|x_j-x_k|\right)\right).
\end{eqnarray*}
\end{proof}
\end{theorem}
%-----------------------------------------------------------------------------------------------------------------------------------------------------------%
Based on Lemma \ref{FK} and inequalities (\ref{FK2}) and (\ref{FK4}), one can employ an argument similar to that used in the proof of Theorem \ref{MT1} to
prove the following related result.
%-----------------------------------------------------------------------------------------------------------------------------------------------------------%
\begin{theorem}\label{MT2}
Let $x_0,\ldots,x_{n-1}$ be $\tau$-measurable operators and
$\alpha_0,\ldots,\alpha_{n-1}$ be positive real numbers such that
$\Sigma_{j=0}^{n-1}\alpha_j=1$. Let $\varphi$ be a continuous
increasing function on $\mathbb{R}_+$ with $\varphi(0)=0$ such that
$\psi(t)=\varphi(\sqrt{t})$ is concave. Then
\begin{eqnarray}\label{ca6}
\sum_{j=0}^{n-1}\alpha_j\tau(\varphi(|x_j|))\leq\tau\left(\varphi\left(\left|\sum_{j=0}^{n-1}\alpha_jx_j\right|\right)\right)+
\sum_{0\leq j<k\leq n-1}\tau\left(\varphi\left(\sqrt{\alpha_j\alpha_k}\left|x_j-x_k\right|\right)\right).
\end{eqnarray}
\end{theorem}
%-----------------------------------------------------------------------------------------------------------------------------------------------------------%
Applying Theorem \ref{MT1} and Theorem \ref{MT2} for the functions
$\varphi(t)=t^p$ $(2\leq p<\infty)$ and $\varphi(t)=t^p$
$(0<p\leq2)$, respectively, we obtain natural generalization of
Clarkson inequalities for $n$-tuples of $\tau$-measurable operators.

%-----------------------------------------------------------------------------------------------------------------------------------------------------------%
\begin{corollary}
Let $x_0,\ldots,x_{n-1}$ be $\tau$-measurable operators and
$\alpha_0,\ldots,\alpha_{n-1}$ be positive real numbers such that
$\Sigma_{j=0}^{n-1}\alpha_j=1$. Then
\begin{eqnarray*}
\left\|\sum_{j=0}^{n-1}\alpha_jx_j\right\|^p_p+\sum_{0\leq j<k\leq n-1}(\alpha_j\alpha_k)^{\frac{p}{2}}\left\|x_j-x_k\right\|^p_p
\leq\sum_{j=0}^{n-1}\alpha_j\left\|x_j\right\|^p_p
\end{eqnarray*}
for $2\leq p<\infty$, and
\begin{eqnarray*}
\sum_{j=0}^{n-1}\alpha_j\left\|x_j\right\|^p_p\leq\left\|\sum_{j=0}^{n-1}\alpha_jx_j\right\|^p_p+\sum_{0\leq j<k\leq n-1}(\alpha_j\alpha_k)^{\frac{p}{2}}\left\|x_j-x_k\right\|^p_p
\end{eqnarray*}
for $0<p\leq2$. In particular
\begin{eqnarray*}
\left\|\sum_{j=0}^{n-1}x_j\right\|^p_p+\sum_{0\leq j<k\leq n-1}\left\|x_j-x_k\right\|^p_p
\leq n^{p-1}\sum_{j=0}^{n-1}\left\|x_j\right\|^p_p
\end{eqnarray*}
for $2\leq p<\infty$, and
\begin{eqnarray*}
n^{p-1}\sum_{j=0}^{n-1}\left\|x_j\right\|^p_p\leq\left\|\sum_{j=0}^{n-1}x_j\right\|^p_p+\sum_{0\leq j<k\leq n-1}\left\|x_j-x_k\right\|^p_p
\end{eqnarray*}
for $0<p\leq2$.
\end{corollary}
%--------------------------------------------------------------------------------------------------------------------------------------------------------------%
Let $\varphi(t)=e^{t^2}-1$. Then $\varphi$ is a continuous
increasing function on $\mathbb{R}_+$, $\varphi(0)=0$ and
$\psi(t)=\varphi(\sqrt{t})=e^t-1$ is convex. Applying Theorem
\ref{MT1} to this spacial function, we get the next result.
%--------------------------------------------------------------------------------------------------------------------------------------------------------------%
\begin{corollary}
Let $x_0,\ldots,x_{n-1}$ be $\tau$-measurable operators and
$\alpha_0,\ldots,\alpha_{n-1}$ be positive real numbers such that
$\Sigma_{j=0}^{n-1}\alpha_j=1$. Then
\begin{eqnarray*}
\tau\left(e^{\left|\sum_{j=0}^{n-1}\alpha_jx_j\right|^2}-\textbf{1}\right)+\sum_{0\leq j<k\leq n-1}\tau\left(e^{\alpha_j\alpha_k|x_j-x_k|^2}-\textbf{1}\right)
\leq\sum_{j=0}^{n-1}\alpha_j\tau\left(e^{|x_j|^2}-\textbf{1}\right).
\end{eqnarray*}
\end{corollary}
%----------------------------------------------------------------------------------------------------------------------------------------------------------------%
Now let $\varphi(t)=\log(t+1)$. Then
$\psi(t)=\varphi(\sqrt{t})=\log(\sqrt{t}+1)$ is concave. Applying
Theorem \ref{MT2} to this function, we obtain the following
corollary.

%----------------------------------------------------------------------------------------------------------------------------------------------------------------%
\begin{corollary}
Let $x_0,\ldots,x_{n-1}$ be $\tau$-measurable operators and
$\alpha_0,\ldots,\alpha_{n-1}$ be positive real numbers such that
$\Sigma_{j=0}^{n-1}\alpha_j=1$. Then
\begin{eqnarray*}
\sum_{j=0}^{n-1}\alpha_j\tau\left(\log(|x_j|+1)\right)\leq\tau\left(\log\left(\left|\sum_{j=0}^{n-1}\alpha_jx_j\right|+1\right)\right)+
\sum_{0\leq j<k\leq
n-1}\tau\left(\log\left(\sqrt{\alpha_j\alpha_k}|x_j-x_k|+1\right)\right)\,.
\end{eqnarray*}
\end{corollary}
%----------------------------------------------------------------------------------------------------------------------------------------------------------------%
\section{Refinement of trace inequalities }

In this section, by using some trace inequalities for $n$-tuples of
$\tau$-measurable operators, we obtain a refinement of the Clarkson
inequalities, which was established by Bahatia and Kittaneh in
\cite{BK1}. To this end, we need the following useful identity
\begin{eqnarray}\label{IBK}
\frac{1}{n}\sum_{k=0}^{n-1}\left|\sum_{j=0}^{n-1}\omega_j^kx_j\right|^2=\sum_{j=0}^{n-1}|x_j|^2,
\end{eqnarray}
where $\omega_0,\ldots,\omega_{n-1}$ are the $n$ roots of unity with $\omega_j=e^{2\pi ij/n}$, $0\leq j\leq n-1$ (see \cite{BK1}).
%---------------------------------------------------------------------------------------------------------------------------------------------------------------%
\begin{theorem}\label{TR1}
Let $x_0,\ldots,x_{n-1}$ be $\tau$-measurable operators and let
$\varphi$ be a continuous increasing function on $\mathbb{R}_+$ with
$\varphi(0)=0$ such that $\psi(t)=\varphi(\sqrt{t})$ is convex. Then
\begin{eqnarray}\label{R1}
\sum_{k=0}^{n-1}\tau\left(\varphi\left(\frac{1}{\sqrt{n}}\left|\sum_{j=0}^{n-1}\omega_j^kx_j\right|\right)\right)
&\leq&\tau\left(\varphi\left(\left(\sum_{j=0}^{n-1}|x_j|^2\right)^{\frac{1}{2}}\right)\right)\\\nonumber
&\leq&\frac{1}{n}\sum_{k=0}^{n-1}\tau\left(\varphi\left(\left|\sum_{j=0}^{n-1}\omega_j^kx_j\right|\right)\right).
\end{eqnarray}
\begin{proof}
Using inequalities (\ref{FK1}), (\ref{FK3}) and identity
(\ref{IBK}), we reach the first inequality in (\ref{R1}).
\begin{eqnarray*}
\sum_{k=0}^{n-1}\tau\left(\varphi\left(\frac{1}{\sqrt{n}}\left|\sum_{j=0}^{n-1}\omega_j^kx_j\right|\right)\right)&=&
\sum_{k=0}^{n-1}\tau\left(\psi\left(\frac{1}{n}\left|\sum_{j=0}^{n-1}\omega_j^kx_j\right|^2\right)\right)\\
&\leq&\tau\left(\psi\left(\frac{1}{n}\sum_{k=0}^{n-1}\left|\sum_{j=0}^{n-1}\omega_j^kx_j\right|^2\right)\right)\\
&=&\tau\left(\psi\left(\sum_{j=0}^{n-1}|x_j|^2\right)\right)\\
&=&\tau\left(\varphi\left(\left(\sum_{j=0}^{n-1}|x_j|^2\right)^{\frac{1}{2}}\right)\right),
\end{eqnarray*}
as desired.

Now using Lemma \ref{FK}, inequality (\ref{FK1}) and identity
(\ref{IBK}) we get the second inequality in (\ref{R1}).
\begin{eqnarray*}
\tau\left(\varphi\left(\left(\sum_{j=0}^{n-1}|x_j|^2\right)^{\frac{1}{2}}\right)\right)&=&
\tau\left(\varphi\left(\left(\frac{1}{n}\sum_{k=0}^{n-1}\left|\sum_{j=0}^{n-1}\omega_j^kx_j\right|^2\right)^{\frac{1}{2}}\right)\right)\\
&=&\tau\left(\psi\left(\frac{1}{n}\sum_{k=0}^{n-1}\left|\sum_{j=0}^{n-1}\omega_j^kx_j\right|^2\right)\right)\\
&\leq&\frac{1}{n}\sum_{k=0}^{n-1}\tau\left(\psi\left(\left|\sum_{j=0}^{n-1}\omega_j^kx_j\right|^2\right)\right)\\
&=&\frac{1}{n}\sum_{k=0}^{n-1}\tau\left(\varphi\left(\left|\sum_{j=0}^{n-1}\omega_j^kx_j\right|\right)\right).
\end{eqnarray*}
\end{proof}
\end{theorem}
%--------------------------------------------------------------------------------------------------------------------------------------------------------------%
Based on inequalities (\ref{FK2}), (\ref{FK4}) and identity
(\ref{IBK}), one can employ an argument similar to that used in the
proof of Theorem \ref{TR1} to get the following related result.
%----------------------------------------------------------------------------------------------------------------------------------------------------------------%
\begin{theorem}\label{TR2}
Let $x_0,\ldots,x_{n-1}$ be $\tau$-measurable operators and let
$\varphi$ be a continuous increasing function on $\mathbb{R}_+$ with
$\varphi(0)=0$ such that $\psi(t)=\varphi(\sqrt{t})$ is concave.
Then
\begin{eqnarray*}\label{R2}
\frac{1}{n}\sum_{k=0}^{n-1}\tau\left(\varphi\left(\left|\sum_{j=0}^{n-1}\omega_j^kx_j\right|\right)\right)\leq
\tau\left(\varphi\left(\left(\sum_{j=0}^{n-1}|x_j|^2\right)^{\frac{1}{2}}\right)\right)\leq\sum_{k=0}^{n-1}\tau\left(\varphi\left(\frac{1}{\sqrt{n}}\left|\sum_{j=0}^{n-1}\omega_j^kx_j\right|\right)\right).
\end{eqnarray*}
\end{theorem}
%---------------------------------------------------------------------------------------------------------------------------------------------------------------%
Applying Theorem \ref{TR1} and Theorem \ref{TR2} for the functions
$\varphi(t)=t^p$ $(2\leq p<\infty)$ and $\varphi(t)=t^p$
$(0<p\leq2)$, respectively, we obtain refinements of the Clarkson
inequalities for $n$-tuples of $\tau$-measurable operators.
%---------------------------------------------------------------------------------------------------------------------------------------------------------------%
\begin{corollary}
Let $x_0,\ldots,x_{n-1}$ be $\tau$-measurable operators. Then
\begin{eqnarray*}
n^{-\frac{p}{2}}\sum_{k=0}^{n-1}\left\|\sum_{j=0}^{n-1}\omega_j^kx_j\right\|_p^p
\leq\left\|\left(\sum_{j=0}^{n-1}|x_j|^2\right)^{\frac{1}{2}}\right\|_p^p
\leq\frac{1}{n}\sum_{k=0}^{n-1}\left\|\sum_{j=0}^{n-1}\omega_j^kx_j\right\|_p^p
\end{eqnarray*}
for $2\leq p<\infty$, and
\begin{eqnarray*}
\frac{1}{n}\sum_{k=0}^{n-1}\left\|\sum_{j=0}^{n-1}\omega_j^kx_j\right\|_p^p
\leq\left\|\left(\sum_{j=0}^{n-1}|x_j|^2\right)^{\frac{1}{2}}\right\|_p^p
\leq n^{-\frac{p}{2}}\sum_{k=0}^{n-1}\left\|\sum_{j=0}^{n-1}\omega_j^kx_j\right\|_p^p
\end{eqnarray*}
for $0<p\leq2$.
\end{corollary}
%----------------------------------------------------------------------------------------------------------------------------------------------------------------%
Applying Theorem \ref{TR1} and Theorem \ref{TR2} to special
functions $\varphi(t)=e^{t^2}-1$ and $\varphi(t)=\log(t+1)$
respectively, we obtain the following results.
\begin{corollary}
Let $x_0,\ldots,x_{n-1}$ be $\tau$-measurable operators. Then
\begin{eqnarray*}
\sum_{k=0}^{n-1}\tau\left(e^{\frac{1}{n}\left|\sum_{j=0}^{n-1}\omega_j^kx_j\right|^2}-\textbf{1}\right)\leq
\tau\left(e^{\sum_{j=0}^{n-1}|x_j|^2}-\textbf{1}\right)\leq
\frac{1}{n}\sum_{k=0}^{n-1}\tau\left(e^{\left|\sum_{j=0}^{n-1}\omega_j^kx_j\right|^2}-\textbf{1}\right).
\end{eqnarray*}
\end{corollary}
%----------------------------------------------------------------------------------------------------------------------------------------------------------------%
\begin{corollary}
Let $x_0,\ldots,x_{n-1}$ be $\tau$-measurable operators. Then
\begin{eqnarray*}
\frac{1}{n}\sum_{k=0}^{n-1}\tau\left(\log\left(\left|\sum_{j=0}^{n-1}\omega_j^kx_j\right|+\textbf{1}\right)\right)&\leq&
\tau\left(\log\left(\left(\sum_{j=0}^{n-1}|x_j|\right)^{\frac{1}{2}}+\textbf{1}\right)\right)\\
&\leq&\frac{1}{n}\sum_{k=0}^{n-1}\tau\left(\log\left(\frac{1}{n}\left|\sum_{j=0}^{n-1}\omega_j^kx_j\right|+\textbf{1}\right)\right).
\end{eqnarray*}
\end{corollary}
%-----------------------------------------------------------------------------------------------------------------------%
\section{A general parallelogram law for $\tau$-measurable operators }
The parallelogram law for $n$--tuples of operators has been recently
investigated by the second author \cite{MO}. In this section, by
using some trace inequalities for $n$-tuples of $\tau$-measurable
operators, we obtain an extension of the main result of \cite{MO}
for $n$--tuples of $\tau$-measurable operators.
%------------------------------------------------------------------------------------------------------------------------------%
Let $x_0,\ldots,x_{n-1}, y_1,\ldots, y_n$ be $\tau$-measurable operators and
$\alpha_0,\ldots,\alpha_{n-1}$ be positive real numbers. Then the following
identity can be given similar to \cite[Corollary 2.3]{MO},
\begin{eqnarray}\label{MO1}
\sum_{0\leq i<j\leq
n-1}\left|\sqrt{\frac{\alpha_i}{\alpha_j}}x_i-\sqrt{\frac{\alpha_j}{\alpha_i}}x_j\right|^2+\sum_{0\leq i<j\leq
n-1}\left|\sqrt{\frac{\alpha_i}{\alpha_j}}y_i-\sqrt{\frac{\alpha_j}{\alpha_i}}y_j\right|^2\\=
\sum_{i,j=0}^{n-1}\left|\sqrt{\frac{\alpha_i}{\alpha_j}}x_i-\sqrt{\frac{\alpha_j}{\alpha_i}}y_j\right|^2-\left|\sum_{i=0}^{n-1}(x_i-y_i)\right|^2.\nonumber
\end{eqnarray}
If we set $y_1=\ldots= y_n=0$ in (\ref{MO1}) and $\Sigma_{j=0}^{n-1}\frac{1}{\alpha_j}=1$,
then the following
identity can be given
\begin{eqnarray}\label{MO2}
\sum_{0\leq i<j\leq
n-1}\left|\sqrt{\frac{\alpha_i}{\alpha_j}}x_i-\sqrt{\frac{\alpha_j}{\alpha_i}}x_j\right|^2=
\sum_{i=0}^{n-1}\alpha_i|x_j|^2-\left|\sum_{i=0}^{n-1}x_i\right|^2.
\end{eqnarray}
%------------------------------------------------------------------------------------------------------------------------------%
\begin{theorem}\label{TL}
Let $x_0,\ldots,x_{n-1}, y_1,\ldots, y_n$ be $\tau$-measurable operators and
$\alpha_0,\ldots,\alpha_{n-1}$ be positive real number such that $\sum_{i,j=0}^{n-1}\frac{1}{\sqrt{\alpha_i\alpha_j}}=1$.  Let $\varphi$ be a
continuous increasing function on $\mathbb{R}_+$ with $\varphi(0)=0$
such that $\psi(t)=\varphi(\sqrt{t})$ is convex. Then
\begin{eqnarray}\label{LA}
\sum_{i,j=0}^{n-1}\frac{1}{\sqrt{\alpha_i\alpha_j}}\tau(\varphi(|\alpha_ix_i-\alpha_jy_j|))
&\geq&\\\nonumber
\sum_{0\leq i<j\leq n-1}\tau\left(\varphi\left(\left|\sqrt{\frac{\alpha_i}{\alpha_j}}x_i-\sqrt{\frac{\alpha_j}{\alpha_i}}x_j\right|\right)\right)
&+&\sum_{0\leq i<j\leq n-1}\tau\left(\varphi\left(\left|\sqrt{\frac{\alpha_i}{\alpha_j}}y_i-\sqrt{\frac{\alpha_j}{\alpha_i}}y_j\right|\right)\right)\\\nonumber
&+&\tau\left(\varphi\left(\left|\sum_{i=0}^{n-1}(x_i-y_i)\right|\right)\right).
\end{eqnarray}
If $\psi(t)$ is concave, then the reverse of inequality (\ref{LA})
holds.
\begin{proof}
\begin{eqnarray*}
&&\hspace{-1cm}\sum_{i,j=0}^{n-1}\frac{1}{\sqrt{\alpha_i\alpha_j}}\tau(\varphi(|\alpha_ix_i-\alpha_jy_j|))\\
&=&\sum_{i,j=0}^{n-1}\frac{1}{\sqrt{\alpha_i\alpha_j}}\tau\left(\psi\left(\left|\alpha_ix_i-\alpha_jy_j\right|^2
\right)\right)\\
&\geq&\tau\left(\psi\left(\sum_{i,j=0}^{n-1}\frac{1}{\sqrt{\alpha_i\alpha_j}}\left|\alpha_ix_i-\alpha_jy_j\right|^2
\right)\right)\qquad\qquad\qquad\qquad\qquad\qquad\qquad\qquad(\textrm{by~} (\ref{FK1}))\\
&=&\tau\left(\psi\left(\sum_{0\leq i<j\leq
n-1}\left|\sqrt{\frac{\alpha_i}{\alpha_j}}x_i-\sqrt{\frac{\alpha_j}{\alpha_i}}x_j\right|^2+\sum_{0\leq i<j\leq
n-1}\left|\sqrt{\frac{\alpha_i}{\alpha_j}}y_i-\sqrt{\frac{\alpha_j}{\alpha_i}}y_j\right|^2\right)
+\left|\sum_{i=0}^n(x_i-y_i)\right|^2\right)\\
&&\qquad\qquad\qquad\qquad\qquad\qquad\qquad\qquad\qquad\qquad\qquad\qquad\qquad\qquad\qquad\qquad(\textrm{by~} (\ref{MO1}))\\
&\geq&\sum_{i,j=0}^{n-1}\tau\left(\psi\left(\left|\sqrt{\frac{\alpha_i}{\alpha_j}}x_i-\sqrt{\frac{\alpha_j}{\alpha_i}}x_j\right|^2\right)\right)
+\sum_{i,j=0}^{n-1}\tau\left(\psi\left(\left|\sqrt{\frac{\alpha_i}{\alpha_j}}y_i-\sqrt{\frac{\alpha_j}{\alpha_i}}y_j\right|^2\right)\right)\\
&\quad+& \tau\left(\psi\left(\left|\sum_{i=0}^n(x_i-y_i)\right|^2\right)\right)
\qquad\qquad\qquad\qquad\qquad\qquad\qquad\qquad\qquad\qquad(\textrm{by~} (\ref{FK3}))\\
&=&\sum_{i,j=0}^{n-1}\tau\left(\varphi\left(\left|\sqrt{\frac{\alpha_i}{\alpha_j}}x_i-\sqrt{\frac{\alpha_j}{\alpha_i}}x_j\right|\right)\right)
+\sum_{i,j=0}^{n-1}\tau\left(\varphi\left(\left|\sqrt{\frac{\alpha_i}{\alpha_j}}y_i-\sqrt{\frac{\alpha_j}{\alpha_i}}y_j\right|\right)\right)\\
&\quad +&\tau\left(\varphi\left(\left|\sum_{i=0}^n(x_i-y_i)\right|\right)\right)
\end{eqnarray*}
\end{proof}
\end{theorem}
%------------------------------------------------------------------------------------------------------------------------------%
\begin{corollary}\label{TL1}
Let $x_0,\ldots,x_{n-1}$ be $\tau$-measurable operators and
$\alpha_0,\ldots,\alpha_{n-1}$ be positive real numbers such that
$\Sigma_{j=0}^{n-1}\frac{1}{\alpha_j}=1$.  Let $\varphi$ be a
continuous increasing function on $\mathbb{R}_+$ with $\varphi(0)=0$
such that $\psi(t)=\varphi(\sqrt{t})$ is convex. Then
\begin{eqnarray}\label{LA1}
\sum_{j=0}^{n-1}\frac{1}{\alpha_j}\tau(\varphi(|\alpha_ix_i|)\geq
\sum_{0\leq i<j\leq
n-1}\tau\left(\varphi\left(\left|\sqrt{\frac{\alpha_i}{\alpha_j}}x_i-\sqrt{\frac{\alpha_j}{\alpha_i}}x_j\right|\right)\right)+
\tau\left(\varphi\left(\left|\sum_{i=0}^{n-1}x_i\right|\right)\right).
\end{eqnarray}
If $\psi(t)$ is concave, then the reverse of inequality \ref{LA1}
holds.
\begin{proof}
It is sufficient that put $y_1=\ldots= y_n=0$ in Theorem \ref{TL}.
\end{proof}
\end{corollary}
%--------------------------------------------------------------------------------------------------------------------------------------------------------------%
Applying Corollary \ref{TL1} for the function $\varphi(t)=t^p (2\leq p<\infty)$ and $\varphi(t)=t^p (0< p\leq2)$ respectively,
we get the following corollary.
%--------------------------------------------------------------------------------------------------------------------------------------------------------------%
\begin{corollary}
Let $x_0,\ldots,x_{n-1}$ be $\tau$-measurable operators and
$\alpha_0,\ldots,\alpha_{n-1}$ be positive real numbers such that
$\Sigma_{j=0}^{n-1}\frac{1}{\alpha_j}=1$. Then
\begin{eqnarray*}
\sum_{i=0}^{n-1}\alpha_i^{p-1}\|x_i\|^p\geq\sum_{0\leq i<j\leq
n-1}\left\|\sqrt{\frac{\alpha_i}{\alpha_j}}x_i-\sqrt{\frac{\alpha_j}{\alpha_i}}x_j\right\|_p^p+
\left\|\sum_{i=0}^{n-1}x_i\right\|_p^p
\end{eqnarray*}
for $2\leq p<\infty$, and
\begin{eqnarray*}
\sum_{0\leq i<j\leq
n-1}\left\|\sqrt{\frac{\alpha_i}{\alpha_j}}x_i-\sqrt{\frac{\alpha_j}{\alpha_i}}x_j\right\|_p^p+
\left\|\sum_{i=0}^{n-1}x_i\right\|_p^p\geq\sum_{i=0}^{n-1}\alpha_i^{p-1}\|x_i\|^p
\end{eqnarray*}
for $0<p\leq2$.
\end{corollary}

\textbf{Acknowledgement.} The first author would like to thank Professor Qing-Wen Wang for his great hospitality during his visit of Shanghai University in 2013.
%----------------------------------------------------------------------------------------------------------------------------------------------------------------%


\begin{thebibliography}{99}

\bibitem {BK1} R. Bahatia and F. Kittaneh,\textit{ Clarkson inequalities with several operators}, Bull. London Math. Soc. \textbf{36}: 820--832, 2004.

\bibitem{D} J. Dixmier, \textit{ Von Neumann Algebras}, North-Holland Mathematical Library, Vol. \textbf{27},
North-Holland, Amesterdam, 1981.

\bibitem {DDP} P. G. Dodds, T. K. Dodds and B. de Pagter,\textit{ Noncommutative Banach function spaces}, Math. Z. \textbf{201}: 47--57, 1989.

\bibitem {DDP3} P. G. Dodds, T.K. Dodds and B. de Pagter,\textit{ Noncommutative K\"{o}the duality}, Trans. Amer. Math. Soc. \textbf{339}: 717--750,1993.

\bibitem{FK} Th. Fack and H. Kosaki, \textit{ Generalized $s$-numbers of $\tau$-measurable operators},
Pacific J. Math. \textbf{123}: 269--300, 1986.

\bibitem {HK1} O. Hirzallah and F. Kittaneh,\textit{ Noncommutative Clarkson inequalities for unitarily invariant norms },
Pacific J. Math. \textbf{202}: 363--369, 2002.

\bibitem {HK2} O. Hirzallah and F. Kittaneh,\textit{ Noncommutative Clarkson inequalities for n-tuples of operators},
Integral equations Operator theory \textbf{60}: 369--379, 2008.

\bibitem{K} E. Kissin, \textit{ On Clarkson--McCarthy inequalities for n-tuples of operators},
Proc. Amer. Math. Soc. \textbf{135}: 2483--2495, 2007.

\bibitem{MO} M. S. Moslehian, \textit{ An operator extension of the parallelogram law and related norm inequalities },
Math. Inequal. Appl. \textbf{14}(2): 717--725, 2011.

\bibitem{MST} M. S. Moslehian, M. Tominaga and K. S. Saito, \textit{Schatten $p$-norm inequalities related to an extended operator
parallelogram law}, Linear Algebra Appl. \textbf{435}(4): 823--829, 2011.

\bibitem{N} E. Nelson, \textit{ Note on noncommutative integration}, J. Funct. Anal. \textbf{15}: 103--116, 1974.

\bibitem{P} B. de Pagter, \textit{ Noncommutative Banach function spaces}, Positivity, 197-227, Trends Math., Birkh\"user, Basel, 2007.

\bibitem{SAD} Gh. Sadeghi, \textit{Non-commutative Orlicz spaces associated to a
modular on $\tau$-measurable operators}, J. Math. Anal. Appl.
\textbf{395}(2): 705--715, 2012.

\bibitem{S} B. Simon, \textit{ Trace Ideals and Their Applications}, Cambridge University press, Cambridge, 1979.

\bibitem{Y} F. J. Yeadon, \textit{ Noncommutative $L^p$-spaces}, Proc. Camb. Phil. Soc. \textbf{77}: 91--102, 1975

\end{thebibliography}
\end{document}